\documentclass{amsart}

\usepackage{amssymb}
\usepackage[all]{xy}

\newtheorem{thm}{Theorem}

\newtheorem{prop}[thm]{Proposition}

\theoremstyle{definition}
\newtheorem{defn}[thm]{Definition}

\newtheorem{say}[thm]{}
\newtheorem{exmp}[thm]{Example}

\newtheorem{rem}[thm]{Remark}          

\newtheorem*{ack}{Acknowledgments}      

\newtheorem{defn-thm}[thm]{Definition--Theorem}  
\newtheorem{defn-lem}[thm]{Definition--Lemma}  

\theoremstyle{remark}

\renewcommand{\c}[0]{{\mathbb C}}  

\renewcommand{\o}[0]{{\mathcal O}} 
\newcommand{\z}[0]{{\mathbb Z}}

\renewcommand{\r}[0]{{\mathbb R}}

\newcommand{\s}[0]{{\mathbb S}}

\newcommand{\p}[0]{{\mathbb P}}

\newcommand{\map}[0]{\dasharrow}
\newcommand{\qtq}[1]{\quad\mbox{#1}\quad}

\newcommand{\pic}[0]{\operatorname{Pic}}

\newcommand{\chow}[0]{\operatorname{Chow}}

\def\loccoh#1.#2.#3.#4.{H^{#1}_{#2}(#3,#4)}

\DeclareMathAlphabet{\mathchanc}{OT1}{pzc}%
                                {m}{it}

\usepackage[all]{xy}\xyoption{dvips}

\begin{document}
\bibliographystyle{amsalpha}


\title[Gromov--Witten--Welschinger invariants]{Examples of vanishing Gromov--Witten--Welschinger invariants}
\author{J\'anos Koll\'ar}

\maketitle

The aim of this note is to give examples of real enumerative
problems without real solutions.
Methods to define and compute
real enumerative invariants in dimensions 2 and 3 were  developed by
\cite{MR2126497} and partially extended to higher dimensions in
 \cite{2013arXiv1309.4079G}. 
The number of real solutions was determined or estimated in many cases;
see \cite{MR2140988, MR2827813, MR3020146, MR3090704, 2013arXiv1312.2921I, 2013arXiv1309.4079G, 2014arXiv1401.1750G}.
Computations using tropical methods  are given in 
\cite{MR2359091, MR2734349, MR2775877, MR3120740, brug-conic}.

The number of rational curves of degree $d$ defined over $\r$
in $\p^2$ that pass through a given real set  of $3d-1$ points
is almost always nonzero; see the extensive tables in
\cite{MR2775877}. By contrast,   the
predicted number 
of rational curves of degree $d$ in $\p^3$, defined over $\r$,
that pass through a given real set  of $2d$ points
is always  0 for $d$ even; see \cite[2.4]{MR2126497} or
\cite[Cor.1.4]{2013arXiv1309.4079G}.
These formulas count  curves with signs, thus
it is not clear that there are no such curves when the predicted number is 0.
For $d=4$  Mikhalkin found configurations 
of 8 points in $\p^3$ 
without any degree 4 real curve of genus 0 
through them; see \cite[2.6]{MR2126497}.

The main aim of this note is to show that,  for any even $d$,
there are certain types of 
configurations  of $2d$ points in general position with no
 rational curves of degree $d$ defined over $\r$
 passing  through them; see Theorem \ref{no.real.Podd.thm}.
We cover all numerical possibilities where some of the $2d$ points are
real. (We have no such examples for $d$ conjugate point pairs for $d>4$.)

The key to the examples is the study of a degenerate situation when
all the points lie on a degree 4 elliptic curve in $\p^3$.
It turns out that, in this case,  all the curves in question lie on 
some quadric surface; these in turn can be understood by
studying the torsion points on the elliptic curve.
The geometry of the situation is summarized in 
Proposition \ref{Cd.overC.prop}.

In the first non-trivial case, when $d=4$, the  point-octet $P_8$  always
lies  on a degree 4 elliptic curve. In Section \ref{sec.2}
we describe  all possibilities. There are 4 topologically distinct types 
if $P_8$ consists of 4 conjugate pairs of complex points
and 11 topologically distinct types 
if $P_8$ consists of 8 real points.
 Already for $d=6$
the method is unlikely to cover all cases.

More examples of enumerative
problems without real solutions are given in Section \ref{sec.4}.
In most cases we study what happens when all the constraints lie on a
quadric hypersurface. 
This is similar to the method used in \cite{2007arXiv0707.4317W}
for surfaces.
Our examples account for all configurations
involving rational curves of degree $\leq 3$ in $\p^3$.
We also get two  infinite series:
lines in $\p^{4n-1}$ that meet 4 linear subspaces of dimension $2n-1$
(Example \ref{lines.in.4n-1}) 
and rational curves of degree $2d+1$ 
 in $\p^3$ that pass through $4d$ points and 
 meet 4 lines
(Example \ref{higher.in.P3}); the latter pointed out by Zinger.

\section{Degree 4 elliptic curves in $\p^3$}

\begin{say}[Eight points in $\p^3$]\label{8.pts.say}
Let $P_8$ be 8 general points in $\p^3$.
The space of quadrics in $\p^3$ is isomorphic to $\p^9$, hence $P_8$
 lies on a unique pencil of quadrics  $|Q_{\lambda}|$
whose base locus is a degree 4 elliptic curve $E_4$.

Let $C_4\subset\p^3$ be any irreducible, degree 4 curve through $P_8$.
Working over an algebraically closed field, 
any point of $C_4\setminus P_8$ is contained in some
$Q=Q_{\lambda} $. Then $C_4$ and $Q$ meet in $\geq 9$ points,
hence $C_4\subset Q$. A smooth quadric is isomorphic to
$\p^1\times \p^1$ and linear equivalence classes are classified
by their bidegree.
(See Definition \ref{marked.quad.say} for our convention for
marked quadrics and  bidegree.)
Since $C_4$ has degree 4, it
 has bidegree either $(2,2)$ or $(1,3)$ (up to changing the marking of $Q$).
In the first case, $C_4\sim 2H_Q$ where $H_Q$ is the hyperplane class.
Except when $C_4=E_4$, this implies that
 $P_8= (C_4\cap E_4)\sim 2H_Q|_{E_4}$
This is not the case for general $P_8$ as can be seen by varying
one of the points of $P_8$ within $E_4$. 

If the quadric $Q$ is singular then
$2C_4\sim 4H_Q$. As before this would give
$2P_8\sim 4H_Q|_{E_4}$, which is again not the case for general $P_8$.
Thus the quadric $Q$ is smooth
if it contains a curve $C_4\neq E_4$
and, for a suitable choice of the marking,   $C_4$  has bidegree
  $(1,3)$ on $Q$. The   marking of $Q$  gives a line bundle 
$L_2:=\pi_1^*\o_{\p^1}(1)|_{E_4}$
of degree 2 on $E_4$. 
Conversely, a degree 2 line bundle 
$L_2$ on $E_4$ determines a quadric surface that is swept out by the lines
$\langle p_1, p_2\rangle$ where $p_1, p_2\in E_4$ are points such that
$\o_{E_4}(p_1+p_2)\cong L_2$.

Set $H_4:=H_Q|_{E_4}$; this is the hyperplane class on $E_4$,
independent of the choice of $Q$.
The other 
 coordinate projection $\pi_2:Q\to \p^1$ corresponds to  $H_4-L_2$ and so
$$
P_8=(E_4\cdot C_4)\sim 3L_2+(H_4-L_2)=H_4+2L_2.
\eqno{(\ref{8.pts.say}.1)}
$$
Thus a marked quadric $Q\in |Q_{\lambda}|$ 
that contains a curve $C_4$  corresponds to
a solution of the equation
$$
2L_2\sim P_8-H_4\qtq{where} L_2\in \pic(E_4).
\eqno{(\ref{8.pts.say}.2)}
$$
Over an algebraically closed field (of characteristic $\neq 2$)
there are 4 solutions $L_2^1,\dots, L_2^4$ of the equation (\ref{8.pts.say}.2);
any two solutions differ by a 2-torsion point of $\pic(E_4)$.
We get 4 different quadric surfaces $Q^1,\dots, Q^4\in |Q_{\lambda}|$.

We claim that on each  quadric $Q^i$ there is a unique curve
$C_4^i$ that passes through $P_8$. To see this note that a
 curve of bidegree  $(1,3)$ in $\p^1\times \p^1$ can be viewed as the graph
of a degree 3 map
$$
\p^1\to \p^1\qtq{given by} (s{:}t)\mapsto \bigl(g_1(s,t) {:}g_1(s,t)\bigr),
$$
where $g_1, g_2$ are cubic forms. Passing through any given point
gives a linear equation on the 8 coefficients of $g_1, g_2$.
Thus passing through 7  points $p_1,\dots, p_7$ of $P_8$  gives 
a unique pair
$g_1, g_2$, up to scalar.

A nice feature is that passing through the 8th point $p_8$ on $E_4$
comes for free.
Indeed, if the resulting curve passes through an 8th point $q\in E_4$, then
$$
p_1+\cdots+p_7+q= (E_4\cdot C_4) \sim 3L_2+(H_4-L_2)\sim 
p_1+\cdots+p_7+p_8
$$
shows that $q\sim p_8$ hence $q=p_8$ since $E_4$ is elliptic.
Thus  get  4 curves $C_4^i\subset Q^i$ for $ i=1,\dots, 4$.

Furthermore,  a curve $C_4^i$ 
is defined over a subfield $k$ $\Leftrightarrow$  the
corresponding quadric surface $Q^i$ is defined over $k$
$\Leftrightarrow$ the
corresponding $[L_2^i]\in \pic(E_4)$  is a $k$-point.
(Note that even for $k=\r$ it can happen that
a $k$-point of $\pic(E_4)$ does not correspond to an actual line bundle
 on $E_4$ that is defined over $k$.
In this case necessarily $E_4(k)=\emptyset$; see (\ref{pic.R.ell}.3).)
\end{say}

For $d>4$, general point sets $P_{2d}\subset \p^3$
do not lie on a degree 4 elliptic curve. However,
it turns out that point sets $P_{2d}$ that do lie 
on a degree 4 elliptic curve can be studied the same way.
The equation (\ref{8.pts.say}.2) is replaced by a similar one.
A new complication that arises is that on a given
quadric surface $Q^i$ there are usually many curves $C_d^{ij}$.

\begin{defn}\label{marked.quad.say}
A {\it marked  quadric surface} $Q\subset \p^3$
is a smooth quadric plus a choice of  a coordinate  projection 
$\pi_1:Q\cong \p^1\times \p^1\to \p^1$.
The other coordinate  projection is denoted by $\pi_2$. 
 This choice is
equivalent to fixing an isomorphism $\pic(Q)\cong \z^2$
such that both $(1,0)$ and $(0,1)$ correspond to lines.

We say that a curve $B$ on $Q$ has type $(a,b)$ if
$\o_Q(B)\cong \pi_1^*\o_{\p^1}(a)\otimes \pi_2^*\o_{\p^1}(b)$.
\end{defn}

\begin{prop} \label{Cd.overC.prop}
Let $k$ be an algebraically closed field 
and $|Q_{\lambda}|$   a pencil of quadrics  in $\p^3$
whose base locus is an elliptic curve $E$.
Let $P_{2d}\subset E$ be a set of $2d$ general points 
and $C_d$  a connected curve  of degree $d$ through $P_{2d}$
not containing $E$.

Then $C_d$ is irreducible and is contained in one of the
quadrics $Q=Q_{\lambda}$.
Furthermore, such marked quadrics $Q=Q(a, L_2)$ 
correspond to
\begin{enumerate}
\item a choice of $0<a<d/2$ and
\item a solution of $(d-2a)L_2\sim P_{2d}-aH_4$ where $L_2\in \pic_2(E)$
and $H_4$ is the hyperplane class on $E$.
\end{enumerate}
For a fixed  marked  quadric $Q=Q(a, L_2)$  there can be many such curves $C_d$
but they all have type $(a,d-a)$ in $\pic(Q)$.
\end{prop}

Proof. Assume first that
 $C_d$ is irreducible.
Any point of $C_d\setminus P_{2d}$ is contained in some
$Q=Q_{\lambda} $. Then $C_d$ and $Q$ meet in $\geq 2d+1$ points,
hence $C_d\subset Q$. 
Thus $C_d$ has type $(a,d-a)$ for some $0<a\leq d/2$
and suitable choice of the coordinate  projection 
$Q\cong \p^1\times \p^1\to \p^1$.

The choice of $Q$ plus 
the projection $Q\to \p^1$ corresponds to a line bundle
$L_2\in \pic_2(E)$. The other 
projection $Q\to \p^1$ corresponds to  $H_4-L_2$ and so
$$
P_{2d}=(E\cdot C_d)\sim (d-a)L_2+a(H_4-L_2)=aH_4+(d-2a)L_2.
\eqno{(\ref{Cd.overC.prop}.3)}
$$
If $a=d/2$ then this gives $P_{2d}\sim aH_4 $, which is not
the case for general $P_{2d}$. 
Thus the case  $a=d/2$ is excluded and the rest of
(1--2) is clear.

In order to complete the proof, we need to exclude the
reducible cases. Let now $\sum C^i$ be a degree $d$,
connected but 
possibly reducible curve passing through $P_{2d}$.
If one of the $C^i$ passes through more than $2\deg C^i$ points
then it has $>2\deg C^i$ intersection with every $Q_{\lambda}$,
thus $C^i=E$. Otherwise, every  $C^i$ 
passes through exactly $2\deg C^i$ points of $P_{2d}$, thus
it is a curve as described  above.
Furthermore, different $C^i$ pass through different points of
$P_{2d}$. 

If two curves $C^i, C^j$ are contained in different 
 quadric surfaces $Q^i\neq Q^j$, then
the curve is disconnected  since
$C^i\cap C^j\subset Q^i\cap Q^j=E_4$ but the different $C^i$
pass through different points. Thus $\sum C^i$ is contained in
a single quadric surface $Q$. 

If $C^i$ is of type  $(a^i, d^i-a^i)$ then we get subsets
$P^i\subset P_{2d}$ of degree $2d^i$ such that
$$
 P^i\sim (d^i-2a^i)L_2 + a^iH_4.
\eqno{(\ref{Cd.overC.prop}.4)}
$$
Together with
(\ref{Cd.overC.prop}.3) we get
$$
(d-2a)P^i\sim (d^i-2a^i)P_{2d}+ \bigl(a^i(d-2a)-a(d^i-2a^i)\bigr)H_4.
\eqno{(\ref{Cd.overC.prop}.5)}
$$
If $k$ is not an algebraic closure of a finite field, then
we can choose the points $p_i$ such that 
$p_1,\dots, p_{2d}$ and $H_4$ are
linearly independent in $\pic(E_4)$, thus
(\ref{Cd.overC.prop}.5) is impossible unless $P^i=P_{2d}$.
Thus there are no reducible curves $C_d$ for a
very general choice of the points $P_{2d}$.
The conclusion then holds in a Zariski open set of
$E^{2d}$. 
(If $k$ is  an algebraic closure of a finite field
then first we find $P_{2d}$ over the  algebraic closure of $k(x)$
and the note that the above Zariski open set of
$E^{2d}$ must contain $k$-points.)
\qed

\begin{rem} The above method reduces the computation of 
GW-invariants  on $\p^3$  (with 0-dimensional constraints) 
to GW-invariants on $\p^1\times \p^1$. 
It would be interesting to work this out in detail.

This should work over the reals as well, but there is a
subtle point that I find confusing. Following the method
gives that, on the quadric surface $Q$ we need 
to find curves of degree $d$ passing through $P_{2d}$.
However, the space of rational curves of degree $d$ on $Q$
has dimension $2d-1$, thus GW-theory counts the number of
curves that pass through a subset $P_{2d-1}\subset P_{2d}$.
As we noted at the end of Paragraph \ref{8.pts.say}, these curves
then automatically pass through the last point of $P_{2d}$ as well.
However, in the  case when we start with
conjugate point pairs, we can not choose $P_{2d-1}$ to be real.
Nonetheless, the answer has a real structure.

I also have not computed  the normal bundle of the resulting curves;
 this also effects the count.
\end{rem}

\section{Real point-octets in $\p^3$}\label{sec.2}

Starting with a real point-octet in $\p^3$, 
we analyze the number and types of real curves $C_4^i$ 
in (\ref{8.pts.real.say}) to get a complete list of possibilities.

\begin{say}[Picard group of a real elliptic curve]\label{pic.R.ell}
Let $E$ be a real elliptic curve. The Picard group 
of the corresponding complex  elliptic curve is denoted by $\pic=\pic(E)$.
Let $\pic_r\subset \pic$ denote the set of
complex line bundles of degree $r$.

The set of real points of $\pic$ is denoted by $\pic(\r)$. 
Corresponding to the 3 different topological types of $E(\r)$,
there are 3 different descriptions for $\pic(\r)$. 

(\ref{pic.R.ell}.1)  $E(\r)\sim \s^1$. 
In this case $\pic_r(\r)\sim \s^1$ for every $r$.

(\ref{pic.R.ell}.2) $E(\r)\sim\s^1\amalg \s^1$.
In this case $\pic_r(\r)\sim\s^1\amalg \s^1$ for every $r$.

If $r$ is even, then one of these components,  denoted by
$\pic_r^0(\r)$ consists of line bundles that are 
topologically trivial on both components of  $E(\r)$.
We call this component and the line bundles in it  {\it even.} 
The other component,  denoted by
$\pic_r^1(\r)$ consists of line bundles that are 
topologically nontrivial on both components of  $E(\r)$.
We call this component and the line bundles in it  {\it odd.} 

For any complex line bundle $L$, the tensor product
$L\otimes \bar L$ is even. 
Restricting to $\pic(\r) $ 
this is the same as multiplication by 2, denoted by $m_2$.
Thus the 
 image of $m_2: \pic(\r)\to \pic(\r)$
is the union of the even components.

(\ref{pic.R.ell}.3) $E(\r)=\emptyset$. 
In this case $\pic_r(\r)=\emptyset$ for odd $r$ and 
 $\pic_r(\r)\sim\s^1\amalg \s^1$ for even $r$.

If $r$ is even, then one of these components,  denoted by
$\pic_r^0(\r)$ consists of real line bundles.
For $r=2$ these correspond to degree 2 maps to $\p^1$.
We call this component and the line bundles in it  {\it even.} 
The other component,  denoted by
$\pic_r^*(\r)$ consists of points that do not correspond to
real line bundles. 
For $r=2$ these correspond to degree 2 maps to 
the ``empty'' conic
$\tilde\p^1:= (x_0^2+x_1^2+x_2^2=0)$.
We call this component   {\it twisted.} 
As before, the image of $m_2: \pic(\r)\to \pic(\r)$
is the union of the even components.
\end{say}

\medskip

\begin{say}[Counting quartic curves through 8 points in $\p^3$]
\label{8.pts.real.say}
Let $P_8\subset \p^3$ be a real set of 8 points in general position, that is,
$P_8$ is made up of real points and complex conjugate point pairs.

If $P_8$ is real then so is $E_4$.  
The answer to our problem of counting real curves
$C_4$ of degree 4 passing through $P_8$ is determined by the 
topological type of $E_4(\r)$ and the positions of
$H_4\in \pic_4(\r)$ and $[P_8]\in \pic_8(\r)$.

There are 3 possibilities for $E_4(\r)$.
\medskip 

(\ref{8.pts.real.say}.1)  $E_4(\r)=\s^1$.  Then necessarily
$H_4\in \pic_4^0(\r)$ and $[P_8]\in \pic_8^0(\r)$.
Thus $P_8-H_4\in \pic_4^0(\r)$ and $2L_2\sim P_8-H_4$
has 2 real solutions. This gives 
2 real curves $C_4^1, C_4^2$.

\medskip 

(\ref{8.pts.real.say}.2)  $E_4(\r)=\s^1\amalg \s^1$. 
There are 4 subcases.

a) $H_4$ is odd $P_8$ is even. Then $ P_8-H_4$ is odd, 
thus $2L_2\sim P_8-H_4$
has no real solutions.
There are no degree 4  rational curves $C_4$ defined over $\r$.

b) $H_4$ is even $P_8$ is odd. Again $ P_8-H_4$ is odd and no real curves.

c)  $H_4$ and  $P_8$ are both  even.  Then $ P_8-H_4$ is even, 4 real curves.

d) $H_4$ and  $P_8$ are both  odd.  Then $ P_8-H_4$ is even, 4 real curves.

(Note that a) and c) can happen if there are no real points in $P_8$
and also if   all points of  $P_8$ are real  while
$P_8$ must have at least 2 real points in cases b) and d).)

(\ref{8.pts.real.say}.3) $E_4(\r)=\emptyset$.  Then $P_8$ has no real points
hence $[P_8]\in \pic_8^0(\r)$ and also $H_4\in \pic_4^0(\r)$.
Thus $ P_8-H_4$ is even and  $2L_2\sim P_8-H_4$
has 4 solutions in $\pic_2(\r)$.

Note that 2 of these solutions correspond to real line bundles,
giving quadrics $Q\cong \p^1\times \p^1$ and
curves $C_4\cong \p^1$. The other 2 solutions do not 
correspond to real line bundles. These give 
``empty'' quadrics $Q\cong (x_0^2+x_1^2+x_2^2+x_3^2=0)$ and
``empty'' curves $C_4\cong (x_0^2+x_1^2+x_2^2=0)$.

Thus we get 2 real curves $C_4\cong \p^1$ and 
2 real curves $C_4\cong \tilde\p^1$.
\end{say}

\begin{rem} There is a twisted form of (\ref{8.pts.real.say}.3)
when   $E_4$ is  inside  $\tilde\p^3$.
Necessarily  $E_4(\r)=\emptyset$.
Then  $H_4$ sits in $\pic_4^*(\r)$
but $P_8$ is in $\pic_8^0(\r)$. 
Thus  $2L_2\sim P_8-H_4$ has no real solutions.
This is, however, not surprising since there are
no even degree rational curves in $\tilde \p^3$.
\end{rem}

\section{Even degree rational curves in $\p^3$}

Of the cases studied in Section \ref{sec.2}, only one
yields a simple answer for $d>4$.

\begin{prop}\label{no.real.Podd.prop}
 Let $E\subset \p^3$ be a real elliptic curve of degree 4
such that  $E(\r)=\s^1\amalg \s^1$ and $H_4=\o_{\p^3}(1)|_E$ is even.
(Equivalently,  both components of $E(\r)$ are
trivial in $\pi_1(\r\p^3)$.)
Let  $P_{4d}\subset E$ be  a general real subset that has 
 an odd number of points on both components of $E(\r)$.
Let $C_{2d}$ be a connected  real curve 
of degree $\leq 2d$ that contains $P_{4d}$.
Then $E\subset C_{2d}$.
\end{prop}

Proof. Assume that $E\not\subset C_{2d}$. By (\ref{Cd.overC.prop}),
$C_{2d}$ is geometrically irreducible and it is contained in a
 quadric $Q(a,L_2)$ as in (\ref{Cd.overC.prop}.2).
Since $C_{2d}$ is real, so is the  quadric $Q(a,L_2)$.
  By (\ref{Cd.overC.prop}.2), such quadrics correspond to
the  solutions of the equation
$$
(2d-2a)L_2\sim P_{4d}-aH_4.
\eqno{(\ref{no.real.Podd.prop}.1)}
$$
Here $P_{4d}$ is odd and $H_4$ is even, thus $P_{4d}-aH_4$ is odd.
As we noted in (\ref{pic.R.ell}.2), $P_{4d}-aH_4$ is
not an even multiple of a real point of $\pic(E)$. \qed
\medskip

\begin{rem} Once we have established
 that such a curve $C_{2d}$ must lie on a quadric surface,
a simple topological argument also shows 
that (\ref{no.real.Podd.prop}.1) has no solutions, not even
for homology classes with $\z/2$-coefficients.

Under the natural homeomorphism  $Q(\r)\sim {\mathbb S}^1\times {\mathbb S}^1$
the homology class of each component of $E_i\subset E(\r)$ is either
$(0,0)$ or $(1,1)$ in $H_1\bigl(Q(\r), \z/2\bigr)$.
Since $C_{2d}$ has even degree, the homology class of  $C_{2d}(\r)$ is 
again either
$(0,0)$ or $(1,1)$. Thus  $E_i\cap C_{2d}(\r)$ is always even.
\end{rem}

\begin{thm}\label{no.real.Podd.thm}
 Let $E\subset \p^3$  and $P_{4d}\subset E$ be  as in
(\ref{no.real.Podd.prop}). Then there is a (semi\-algebraic) open subset
$[P_{4d}]\in U\subset (S^{4d}\p^3)(\r)$ such that 
if $R_{4d}$ is a real set of $4d$ points in $\p^3$ and
$[R_{4d}]\in U$, then there is no
connected real curve of degree $\leq 2d$ with geometrically rational 
irreducible components
that contains $R_{4d}$. 
\end{thm}

Proof.  Note that $4d$-element subsets of $\p^3$ are parametrized by
the points of the symmetric power $S^{4d}\p^3 $ and real
$4d$-element subsets correspond to the real points 
of the symmetric power $(S^{4d}\p^3)(\r)$
(which is not the same as the symmetric power 
of the real points   $S^{4d}(\r\p^3)$).

Let $W\subset (S^{4d}\p^3)(\r)$ denote the set of points
$[R_{4d}]$ such that there is a
connected real curve of degree $\leq 2d$ with geometrically rational components
that contains $R_{4d}$.  By (\ref{no.real.Podd.prop}) we know that
$[P_{4d}]\not\in W$. 
We claim that  $W$ is a closed semialgebraic subset of 
$(S^{4d}\p^3)(\r)$. If this holds then we can take
$U:=(S^{4d}\p^3)(\r)\setminus W$.

The claim about $W$ follows a standard argument.
Let $\chow_e(\p^3)$ denote the Chow variety of  curves 
(or 1-cycles) of
degree $\leq e$ in $\p^3$ \cite[Sec.I.3]{rc-book}. Let 
$\operatorname{RatCycles}_e\subset \chow_e(\p^3)$ denote the
subset corresponding to curves that are connected with rational 
irreducible components.  $\operatorname{RatCycles}_e$ is a Zariski closed
subset by \cite[II.2.2]{rc-book}.

Next consider $(S^{4d}\p^3)\times \operatorname{RatCycles}_{2d}\times \p^3$
with coordinate projections $\pi_i$.
Let ${\mathcal U}_{4d}\subset (S^{4d}\p^3)\times\p^3$ be the universal family of
$4d$-element subsets of $\p^3$ and
${\mathcal C}_{2d}\subset\operatorname{RatCycles}_{2d}\times \p^3$
the universal family of 1-cycles  \cite[I.3.21]{rc-book}. 

Let ${\mathcal Y}_{2d}\subset (S^{4d}\p^3)\times \operatorname{RatCycles}_{2d}$
be the set of pairs  $\bigl([R_{4d}], C_{2d}\bigr)$ satisfying 
$R_{4d}\subset C_{2d}$.
Then ${\mathcal Y}_{2d}$ is   a Zariski closed
subset since it is the complement of the Zariski open subset
$$
(\pi_1\times \pi_2)\Bigl(
(\pi_1\times \pi_3)^{-1}\bigl({\mathcal U}_{2d}\bigr)
\setminus (\pi_2\times \pi_3)^{-1}\bigl({\mathcal C}_{4d}\bigr)
\Bigr).
$$
Finally 
$W=\pi_1\bigl( {\mathcal Y}_{2d}(\r)\bigr)\subset (S^{4d}\p^3)(\r)$
is a closed semialgebraic subset. \qed




\section{Other examples with linear constraints}\label{sec.4}

We get  more examples without real solutions by
studying what happens when we choose the 
linear constraints to lie on a quadric hypersurface.

\begin{say}[Linear subspaces on quadric hypersurfaces]\label{2.fams.on.Q}
\cite[Book IV, Sec.XIII.4]{hodge-ped}

Let $Q^{2n}\subset \p^{2n+1}$ be a 
smooth quadric hypersurface over $\c$.
It contains 2 families of $n$-dimensional linear subspaces
and two subspaces $L^n_1, L^n_2$ belong to the same family
iff 
$$
\dim\bigl( L^n_1\cap L^n_2\bigr)\equiv n\mod 2
$$
(where the empty set has dimension $-1$).
 If $n$ is odd then two general linear subspaces
in the same family are disjoint from each other.



We are especially interested in the ``empty'' real 
 quadric  
$$
Q_E^{2n}:=(x_0^2+\cdots+x_{2n+1}^2=0)\subset \p^{2n+1},
\eqno{(\ref{2.fams.on.Q}.1)}
$$
which contains the conjugate pair of $n$-dimensional linear subspaces
$$
L_{\pm}:=(x_0\pm \sqrt{-1}x_1=\cdots = x_{2n}\pm \sqrt{-1}x_{2n+1}=0)
$$
which are disjoint from each other.
Thus, if $n$ is odd then they are members of the same family.
Therefore both families are defined over $\r$. 
If $n$ is even, then $L_{\pm} $ are members of different families,
hence the two families are conjugate.
\end{say}

\begin{exmp} [Lines in $\p^{4n-1}$]\label{lines.in.4n-1}
For every $n\geq 1$ there are generic configurations 
$$
L^{2n-1}_1, \bar L^{2n-1}_1,L^{2n-1}_2, \bar L^{2n-1}_2\subset \p^{4n-1}
$$
such that no real line intersects all 4 subspaces.
\medskip

Proof. Let $ Q_E^{4n-2}\subset \p^{4n-1}$ be the empty quadric
(\ref{2.fams.on.Q}.1).
Take first a special configuration where 
$L^{2n-1}_1, \bar L^{2n-1}_1,L^{2n-1}_2, \bar L^{2n-1}_2\subset Q_E^{4n-2}$
are disjoint members of the  same family.

We claim that there is no real line $L$ that intersects
all 4 subspaces. Indeed, any line that  intersects
all 4 subspaces has 4 points in common with the  quadric  $Q_E^{4n-2}$,
thus it is contained in $Q_E^{4n-2}$. However, $Q_E^{4n-2}$ has no real points
hence the line can not be real.

Since a limit if real lines is a real line, any small perturbation
of the configuration has the required property. \qed
\end{exmp}

\begin{rem}[Lines in $\p^{4n+1}$] 
It has been known classically that in general there are
$r+1$ complex lines intersecting 
4 linear subspaces $L_1^r,\dots, L_4^r\subset \p^{2r+1}$.
Thus if $r=2n$ is even and $L_1^r,\dots, L_4^r$ is a real
configuration then there is 
 at least 1 real line meeting all 4 subspaces.
\end{rem}

\begin{exmp}[Lines in $\p^3$] Another interesting degenerate situation
is given by taking the lines to lie on a  cubic surface in $\p^3$.

Take a real set of 4 points $p_1,\dots, p_4$ in $\p^2$
plus two more points $p_5, p_6$. By blowing up the 6 points,
we get a cubic surface $S_3\subset \p^3$; the first 4 points
give 4 lines $L_1,\dots, L_4\subset S_3$.  Any line meeting these 4 meets the
cubic in 4 points, thus it is contained in it. 

Thus  the 2 lines meeting $L_1,\dots, L_4$
are obtained as the birational transforms of the conics
through the points  $p_1,\dots, p_4, p_5$
resp.\ $p_1,\dots, p_4, p_6$.
If  $p_5, p_6$ are real, we get real lines. If they are complex conjugates,
we get complex conjugate lines.

Note that, unlike in (\ref{lines.in.4n-1}), here
we can choose 4 or 2 of the lines to be real.
\end{exmp}

\begin{exmp}[Conics in $\p^{4n-1}$] \label{conics.in.4n-1}
For every $n\geq 1$ there are generic configurations 
$$
L^{2n-1}_1, \bar L^{2n-1}_1,L^{2n-1}_2, \bar L^{2n-1}_2, 
L^{2n-1}_3, \bar L^{2n-1}_3
\subset \p^{4n-1}
$$
such that no real conic with real points intersects all 6 subspaces.
 (Note that these constraints define a 2-dimensional
moduli space.)
\medskip

Proof.
Take first a special configuration where all 6 subspaces
are disjoint members of the  same family on $ Q_E^{4n-2}\subset \p^{4n-1}$.

A conic $C$ that  intersects
all 6 subspaces has 6 points in common with the  quadric  $Q_E^{4n-2}$.
If $C$ is irreducible (over $\r$) then
 it is contained in $Q_E^{4n-2}$. Since $Q_E^{4n-2}$ has no real points,
$C$ is an empty conic.

If $C$ is reducible (over $\r$) then its irreducible
components are real lines and at least one of them
must be contained in $Q_E^{4n-2}$, which is impossible.

Since a limit of real conics with real points is a real conic with real points,
 any small perturbation
of the configuration has the required property. \qed
\end{exmp}

We can add either a pair  of real subspaces $L^{4n-3}_4, L^{4n-3}_5 $
or a single real subspace $L^{4n-4}$ to the constraints in 
(\ref{conics.in.4n-1}) to get a vanishing GWW-invariant, but
the example does not show what happens if we add a
conjugate pair of  subspaces $L^{4n-3}_4, \bar L^{4n-3}_4 $.
In $\p^3$ a different  example  excludes all real  conics.

\begin{exmp}[Conics in $\p^3$] 
The space of  conics in $\p^3$ has 
dimension 8. Thus, working with conjugate pairs of linear constraints,
we get a GWW-invariant in the following cases
\begin{enumerate}
\item $p_1, \bar p_1, p_2, \bar p_2$,
\item $p_1, \bar p_1, L_1, \bar L_1, L_2, \bar L_2$,
\item $ L_1, \bar L_1, \dots, L_4, \bar L_4$,
\end{enumerate}
Every conic lies in a unique plane and 4 general points of $\p^3$
do not lie in a plane. Thus there are no conics through 4
general pints.

In the remaining 2 cases there are always complex conics, but 
we claim that   there are generic configurations 
such that no real conic 
intersects all of the constraints.

\medskip

{\it Case 2:} $p_1, \bar p_1, L_1, \bar L_1, L_2, \bar L_2\subset \p^{3}$.
\medskip

 We start with the construction over $\c$.
Two points $p_1, p_2$ determine a line and projecting from it gives
$\pi_1:\p^3\map \p^1$. Given two  points $p_1, p_2$ and
two lines $L_1, L_2$, there is a 1-dimensional 
 family of quadrics passing through them; this gives
$\pi_2:\p^3\map \p^1$. The product of these gives a map
$$
\pi:=\pi_1\times \pi_2:\p^3\map \p^1\times \p^1.
$$
The fibers of $\pi$ are the conics that pass through $p_1, p_2$
and intersect $L_1, L_2$.

Given any other line $L_3$, it intersects the fibers of
$\pi_1$ (which are planes) in 1 point and the 
fibers of
$\pi_2$ (which are quadrics) in 2 points. Thus
$\pi(L)\subset \p^1\times \p^1$ is a curve of bidegree $(2,1)$. 
Two curves of bidegree $(2,1)$ intersect in 4 points, giving 4 conics
that pass through 2 points and intersect 4 lines.

If $p_1\cup p_2$ and $L_1\cup L_2$ are real, then 
$\pi$ is defined over $\r$. In $\p^1\times \p^1$ it is easy to
write down examples of two curves of bidegree $(2,1)$
(real or conjugate pairs) that have no real intersections points.
\medskip

{\it Case 3:} $L_1, \bar L_1,\dots, L_4, \bar L_4\subset \p^{3}$.
\medskip

 Again we start with the construction over $\c$.
A line $\ell$ in $\p^3$ determines a projection
$\pi(\ell):\p^3\map \p^1$. If $C\subset \p^3$ is a smooth conic 
 then $\pi(\ell)|_C:C\to \p^1$ has degree $2$
if $C$ is disjoint from $\ell$ and degree 1 or 0 if 
$C$ intersects $\ell$.

Using this for a pair of disjoint lines $\ell_1, \ell_2\subset \p^3$
we get a map
$$
\pi:=\pi(\ell_1)\times\pi(\ell_2):\p^3\map \p^1\times \p^1.
$$
The fibers of $\pi$ are lines that connect a point of $\ell_1$ to 
a point of $\ell_2$. 
Thinking of $\p^1\times \p^1 $ as a quadric surface in $\p^3$,
the resulting map $\pi:\p^3\map \p^3$ is given by the
quadratic forms on $\p^3$ that vanish on both lines.

Choose 8 lines $L_1, \dots,  L_8 $
to be fibers of $\pi$ in general position.
Then a conic $C$ that intersects all 8 lines corresponds to a
rational curve of bidegree $(2,2)$ on $ \p^1\times \p^1 $ 
passing through 8 general points. However, the space of
rational curves of bidegree $(2,2)$ has dimension 7, thus
there are no such curves through 8 general points.
Thus we conclude that if a degree 2 curve $C\subset \p^3$ meets all 8 lines
then either $\ell_1$ or $\ell_2$ is an irreducible component of $C$.

To get an example over $\r$, use the above construction
 for a conjugate pair of disjoint lines $\ell,\bar \ell\subset \p^3$.
Choose linear forms $\alpha, \beta$ such that 
$\ell=(\alpha= \beta =0)$. Then $\bar \ell=(\bar \alpha= \bar \beta =0)$
and the space of quadratic forms that vanish on both  $\ell,\bar \ell$
is spanned by
$\alpha\bar \alpha, \beta\bar \beta, 
\alpha\bar \beta,  \beta\bar \alpha$.
These satisfy the obvious equation
$$
(\alpha\bar \alpha)(\beta\bar \beta)=(\alpha\bar \beta)(\beta\bar \alpha).
$$
To get a real basis, we change to
$$
\langle\alpha\bar \alpha+ \beta\bar \beta,
\alpha\bar \alpha- \beta\bar \beta,
\alpha\bar \beta+ \beta\bar \alpha, 
\sqrt{-1}(\alpha\bar \beta- \beta\bar \alpha).
$$
These satisfy the equation
$$
(\alpha\bar \alpha+ \beta\bar \beta)^2-
(\alpha\bar \alpha- \beta\bar \beta)^2
=
(\alpha\bar \beta+ \beta\bar \alpha)^2+
\bigl(\sqrt{-1}(\alpha\bar \beta- \beta\bar \alpha)\bigr)^2.
$$
Thus $\pi$ can be thought of as a map
$$
\pi:=\pi(\ell)\times\pi(\bar \ell):\p^3\map  Q,
$$
where $Q$ is isomorphic to the ``sphere'' $(x^2+y^2+z^2=t^2)$.

Now choose $L_1, \bar L_1,\dots, L_4, \bar L_4$
to be fibers of $\pi$ in general position.
By the above considerations,  a degree 2 curve $C$ 
meets all 8 lines iff either $\ell$ or $\bar\ell$ is an irreducible component 
of $C$. 
The only real degree 2 curve 
with this property is $\ell+\bar \ell$. This is, however, 
geometrically disconnected and not a limit of  conics.
\end{exmp}

\begin{exmp}[Cubics in $\p^3$] \label{cubics.in.P3}
The space of rational cubics in $\p^3$ has 
dimension 12. Thus, working with conjugate pairs of linear constraints,
we get a GWW-invariant in the following cases
\begin{enumerate}
\item $p_1, \bar p_1, p_2, \bar p_2,p_3, \bar p_3$,
\item $p_1, \bar p_1, p_2, \bar p_2, L_1, \bar L_1, L_2, \bar L_2$,
\item $p_1, \bar p_1,  L_1, \bar L_1, \dots, L_4, \bar L_4$,
\item $  L_1, \bar L_1, \dots, L_6, \bar L_6$.
\end{enumerate}
It has been classically known that
 there is a unique rational normal curve
through 6 general points in $\p^3$. This curve is real whenever
the 6 points form a real set.

We claim that in the remaining cases there are generic configurations 
such that no degree 3 rational curve defined over $\r$
intersects all of the constraints.

Again we work on the empty quadric
$Q^2_E:=(x_0^2+\cdots+x_3^2=0)\subset \p^3$
and choose special configurations as follows.
\begin{enumerate}
\item[(2')] 
We choose $L_1, \bar L_1, L_2, \bar L_2 $ to be
disjoint members in one family of lines  and the points
in general position on $Q^2_E$.
\item[(3')] 
We choose $L_1, \bar L_1, L_2, \bar L_2, L_3, \bar L_3 $ to be
disjoint members in one family of lines, $ L_4, \bar L_4 $
 are chosen from the other family and the points
in general position on $Q^2_E$.
\item[(4')] 
We choose $L_1, \bar L_1, \dots, L_4, \bar L_4$ to be
disjoint members in one family of lines  and the remaining lines
 to be
disjoint members of  the other family. 
\end{enumerate}

In all of these cases, the following 2 properties hold
\begin{enumerate}
\item[(a)]  there are at least 8 disjoint constraints and
\item[(b)] for both coordinate projections $\pi_i:Q^2_E\to Q^1_E$
at least 4 fibers contain a constraint.
\end{enumerate}
Let $B$ be a real cubic curve that meets all the constraints.
By (a), $B$ and $Q^2_E$ have at least 8 points in common,
thus at least 1 of the irreducible components of $B$
is contained in $Q^2_E$. Since $Q^2_E$ does not contain
odd degree real curves,  $B$
 decomposes as $C+L$ where $C$ is a degree 2
curve contained in $Q^2_E$ and $L$ is a real line.

There are 2 possibilities for $C$.
\begin{enumerate}
\item[(i)] $C$ is a smooth conic.
Note that  $C+L$
can not be written as the image of a geometrically connected 
real curve of arithmetic genus 0 since
we would need to resolve 1 of the nodes, but they
form a conjugate pair. Thus $C+L$
can not be obtained as a limit of real cubics 
of geometric genus 0.
\item[(ii)] $C$ is a  conjugate pair of lines. Then
$Q^2_E\cap (C+L)=C$ but (b) shows that $C$ can not meet all the constraints.
\end{enumerate}

Thus, after a general perturbation we get no real curves. \qed
\end{exmp}

The following generalization of (\ref{cubics.in.P3}.2)
was pointed out by Zinger.

\begin{exmp}[Odd degree rational curves in $\p^3$] \label{higher.in.P3}
For every odd $d\geq 1$ there are generic configurations
$$
p_1, \bar p_1, \dots, p_{d-1}, \bar p_{d-1}, L_1, \bar L_1, L_2, \bar L_2
\subset \p^3
$$
such that 
 no degree $d$ 
rational curve defined over $\r$
intersects all of the constraints. 

Set $P_{2d-2}:=\{p_1, \bar p_1, \dots, p_{d-1}, \bar p_{d-1}\}$. 
As in (\ref{cubics.in.P3}.2') 
we choose $L_1, \bar L_1, L_2, \bar L_2 \subset Q^2_E$ to be
disjoint members in one family of lines  and the points
in general position on $Q^2_E$. 
Let $C_d$ be a 
real curve of degree $d$
with geometrically rational irreducible components
that intersects all of the constraints. 

Note that $C_d$ and $Q^2_E$ have at least $(2d-2)+4$ points in common,
thus at least 1 irreducible component
 of $C_d$ is contained in  $Q^2_E$.
Every real curve contained in $Q^2_E$ has even degree,  thus $C_d$
can not be contained in $Q^2_E$.
We can thus write $C_d=C_{2e}+C_{d-2e}$ where  
$C_{2e}$ is contained in $Q^2_E$ and
none of  the irreducible components of $C_{d-2e}$
is contained in $Q^2_E$. The subscripts indicate the degree.

For general choice of the points, $C_{2e}$
can pass through at most a $4e-2$  element subset  $P_{4e-2}$ of $P_{2d-2}$
and $C_{d-2e}$  can pass through at most a $2d-4e$  element
 subset  $P_{2d-4e}$ of $P_{2d-2}$. Thus $P_{2d-2}$ is a
disjoint union  $P_{4e-2}\cup P_{2d-4e}$. Note further that
$C_{d-2e}\cap Q^2_E=P_{2d-4e}$ which is disjoint from $C_{2e}$.
Thus $C_d=C_{2e}+C_{d-2e}$ is disconnected and it
is not a limit of  geometrically connected  real curves.\qed
\end{exmp}

\cite{2007arXiv0707.4317W} shows that the GWW-invariants for curves 
of any degree in $\p^2$ give the optimal value if 
the constraints (with at most 1 exception) lie near an empty conic. 
The results of \cite{2007arXiv0707.4317W} apply to many
other surfaces as well.
The next examples  illustrate this approach by 
considering the two known vanishing
GWW-invariants for curves in $\p^2$; see  the tables in
\cite{MR2775877}.

\begin{exmp}[Cubics in $\p^2$] \label{deg3.in.Pn.exmp}
There are generic configurations 
$p_1, \bar p_1,\dots , p_4, \bar p_4\in \p^2$
such that no degree 3 rational curve defined over $\r$
passes through all 8 points.
\medskip

Proof. First choose all 8 points on the empty conic
$Q^1_E:=(x^2+y^2+z^2=0)$. Any cubic that
contains the 8 points is of the form $Q^1_E+L$ where $L$
is a line. 
This leads to a contradiction as in (\ref{cubics.in.P3}).\qed
\end{exmp}

\begin{exmp}[Quartics in $\p^2$]\label{deg4.in.P2.exmp}
 There are generic configurations 
$p_1, \bar p_1,\dots , p_5, \bar p_5\in \p^2$
such that there is 
no degree 4 map   $\p^1\to \p^2$ defined over $\r$ whose image 
passes through all 10 points.
 (Note that these constraints define a 1-dimensional
moduli space.)
\medskip

Proof.  First choose all 10 points on the empty conic
$Q^1_E:=(x^2+y^2+z^2=0)$. Then any quartic that
contains the 10 points is of the form $Q^1_E+Q'$ where $Q'$
is a conic. 

A quick case analysis shows that  $Q^1_E+Q'$ 
 can  be written as the image of a real curve of arithmetic genus 0
only when $Q'$ is a conjugate pair of lines $L+\bar L$.
(In this case we can remove the singular point of $Q'$
and one each of the intersections $Q^1_E\cap L, Q^1_E\cap \bar L$.)

Thus $Q^1_E+Q'$ can only be obtained as a limit of real quartics 
of geometric genus 0 if their normalization has no real points.\qed
\end{exmp}

\begin{rem}  A  degeneration argument as in  (\ref{deg4.in.P2.exmp})
 fails to work for higher degree curves.
The curve  $2Q^1_E+L$ is the image of a genus 0
curve consisting of $L$ and a conjugate pair of complex conics, both
mapping isomorphically to $Q^1_E$.
\end{rem}

 \begin{ack}
I thank E.~Brugall\'e, P.~Georgieva, V.~Kharlamov, G.~Mikhalkin and A.~Zinger
 for  answering  my questions and
suggesting many improvements. 
Partial financial support   was provided  by  the NSF under grant number 
DMS-0968337.
\end{ack}


\def\cprime{$'$} \def\cprime{$'$} \def\cprime{$'$} \def\cprime{$'$}
  \def\cprime{$'$} \def\cprime{$'$} \def\dbar{\leavevmode\hbox to
  0pt{\hskip.2ex \accent"16\hss}d} \def\cprime{$'$} \def\cprime{$'$}
  \def\polhk#1{\setbox0=\hbox{#1}{\ooalign{\hidewidth
  \lower1.5ex\hbox{`}\hidewidth\crcr\unhbox0}}} \def\cprime{$'$}
  \def\cprime{$'$} \def\cprime{$'$} \def\cprime{$'$}
  \def\polhk#1{\setbox0=\hbox{#1}{\ooalign{\hidewidth
  \lower1.5ex\hbox{`}\hidewidth\crcr\unhbox0}}} \def\cdprime{$''$}
  \def\cprime{$'$} \def\cprime{$'$} \def\cprime{$'$} \def\cprime{$'$}
\providecommand{\bysame}{\leavevmode\hbox to3em{\hrulefill}\thinspace}
\providecommand{\MR}{\relax\ifhmode\unskip\space\fi MR }
\providecommand{\MRhref}[2]{%
  \href{http://www.ams.org/mathscinet-getitem?mr=#1}{#2}
}
\providecommand{\href}[2]{#2}

\noindent Princeton University, Princeton NJ 08544-1000

{\begin{verbatim}kollar@math.princeton.edu\end{verbatim}}

\end{document}